\documentclass[oneside,english]{amsart}
\usepackage[T1]{fontenc}
\usepackage[latin9]{inputenc}
\usepackage{amsthm}
\usepackage{amssymb}
\usepackage[pdftex]{graphicx}

\makeatletter
\numberwithin{equation}{section}
\numberwithin{figure}{section}

\makeatother

\usepackage{babel}
\begin{document}

\title[\tiny{non-degeneracy of harmonic structures on Sierpinski Gaskets}]{A topological proof of the non-degeneracy of harmonic structures on Sierpinski Gaskets}
%    author one information
\author{Shiping Cao}
\address{Department of Mathematics, Cornell University, Ithaca, 14853, U.S.A.}
\curraddr{} \email{sc2873@cornell.edu}
\thanks{}

%    author two information
\author{Hua Qiu}
\address{Department of Mathematics, Nanjing University, Nanjing, 210093, P. R. China.}
\curraddr{} \email{huaqiu@nju.edu.cn}
\thanks{The research of the second author was supported by the Nature Science Foundation of China, Grant 11471157.}

%    For articles to be published after 1 January 2010, you may use
%    the following version:
%\subjclass[2010]{Primary }


\date{}

\dedicatory{}
\begin{abstract} We present a direct proof of the non-degeneracy of the harmonic structures on the level-$n$ Sierpinski gaskets for any $n\geq 2$, which was conjectured by Hino in [H1,H2] and confirmed to be true by Tsougkas [T] very recently using Tutte's spring theorem. 

\end{abstract}
\maketitle

\section{Introduction }
The theory of analysis on fractals, analogous to that on manifolds, has been being well developed. The pioneering work, developed by Kigami [K1,K2], is the analytic construction of the Laplacians, for a class of finitely ramified fractals, named \textit{p.c.f. self-similar sets}, including the \textit{Sierpinski gasket} as a typical example,  where Laplacians are defined as renormalized limits of graph Laplacians, playing the role of differential operators of second order on manifolds. 

The harmonic functions on fractal domains may have some different nature from the classical ones as a consequence of the finitely ramified topology of the fractals. For example,  the Hexagasket and the Vicsek set [S2]  consume \textit{degenerate harmonic structures} so that  nonconstant harmonic functions vanishing locally on small cells exist.  It seems that such phenomenon always happens on those fractals containing \textit{nonjunction} inner vertices.

Recently, Tsougkas [T] gives a proof on the non-degeneracy of harmonic structures on the \textit{level-$n$ Sierpinski Gasket $\mathcal{SG}_n$}, $n\geq 2$, based on certain geometric graph theory, in which Tutte's spring theorem plays the key role, which was conjectured and numerically checked  for the standard case  for $n\leq 50$ by Hino [H1,H2].

 \textbf{Theorem 1.1.} \textit{For $n\geq 2$, each harmonic structure of $\mathcal{SG}_n$ is non-degenerate.}
 
 In this note, we aim to provide an elementary proof for this theorem using only the basic facts of the discrete Laplacians on finite sets.

For convenience, we list some basic concepts and facts on the \textit{Laplacians on finite sets} below. Readers can find  any unexplained details in the book [K3].

\textbf{Definition 1.2.} \textit{Let $V$ be a finite set. A symmetric linear operator (matrix) $H:l(V)\to l(V)$ is called a Laplacian on $V$ if it satisfies} 

\textit{(1) $H$ is non-positive definite,}

\textit{(2) $Hu=0$ if and only if $u$ is a constant on $V$,}

\textit{(3) $H_{pq}\geq 0$ for all $p\neq q$.}

\textit{Denote $\mathcal{LA}(V)$ the set of Laplacians on $V$.}

\begin{figure}[h]
\begin{center}
\includegraphics[width=3.7cm]{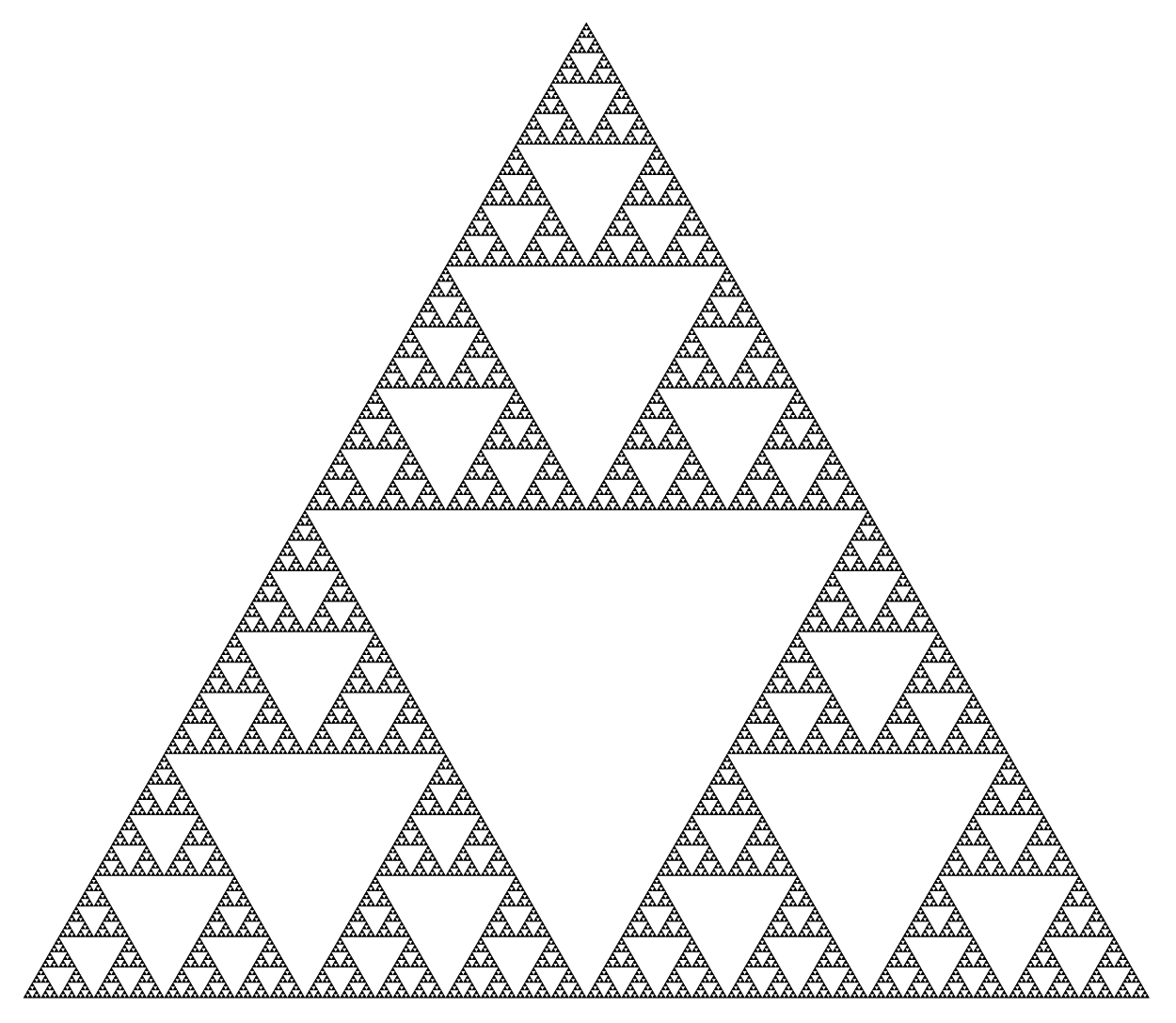}\hspace{0.5cm}
\includegraphics[width=3.7cm]{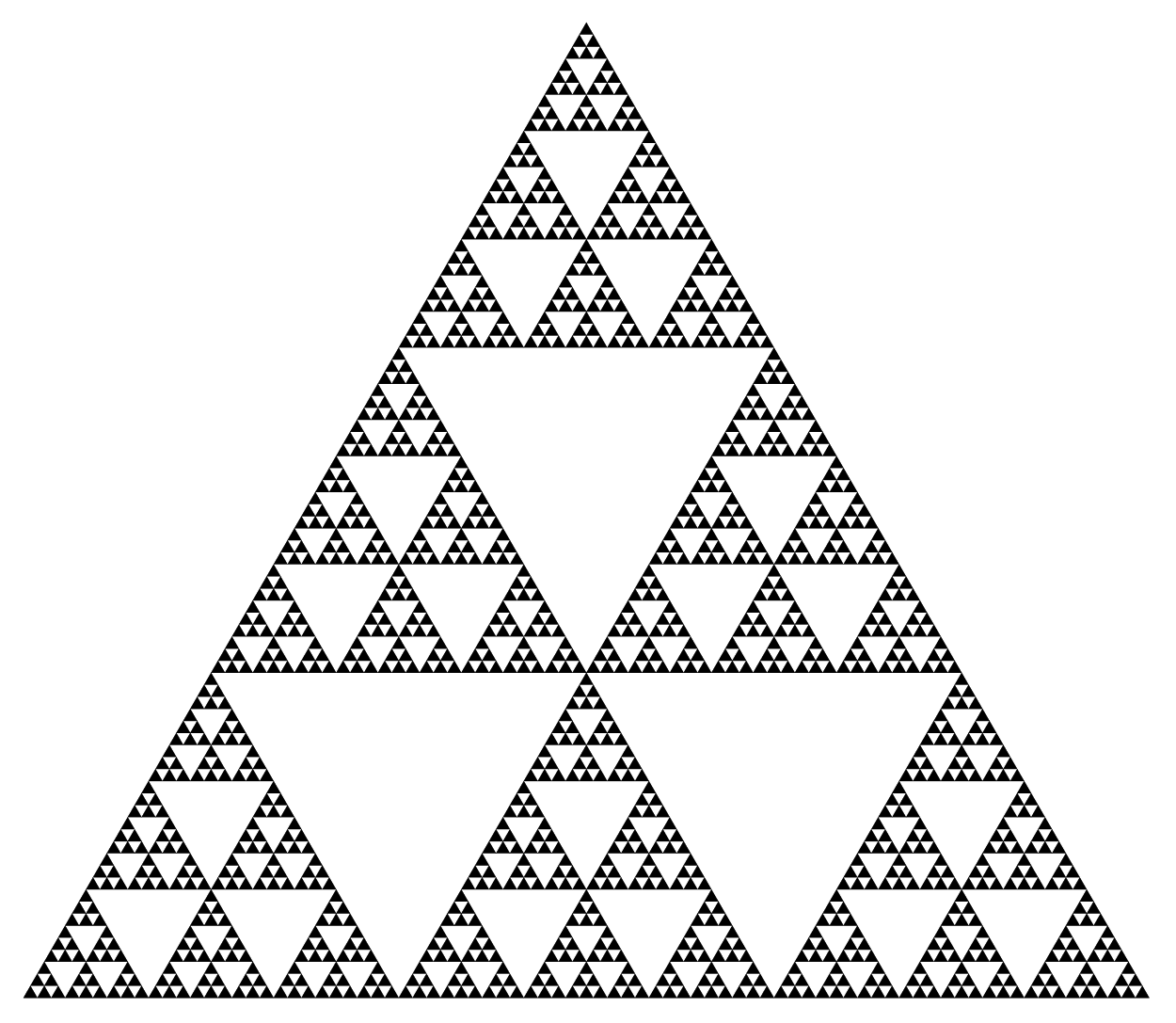}\hspace{0.5cm}
\includegraphics[width=3.7cm]{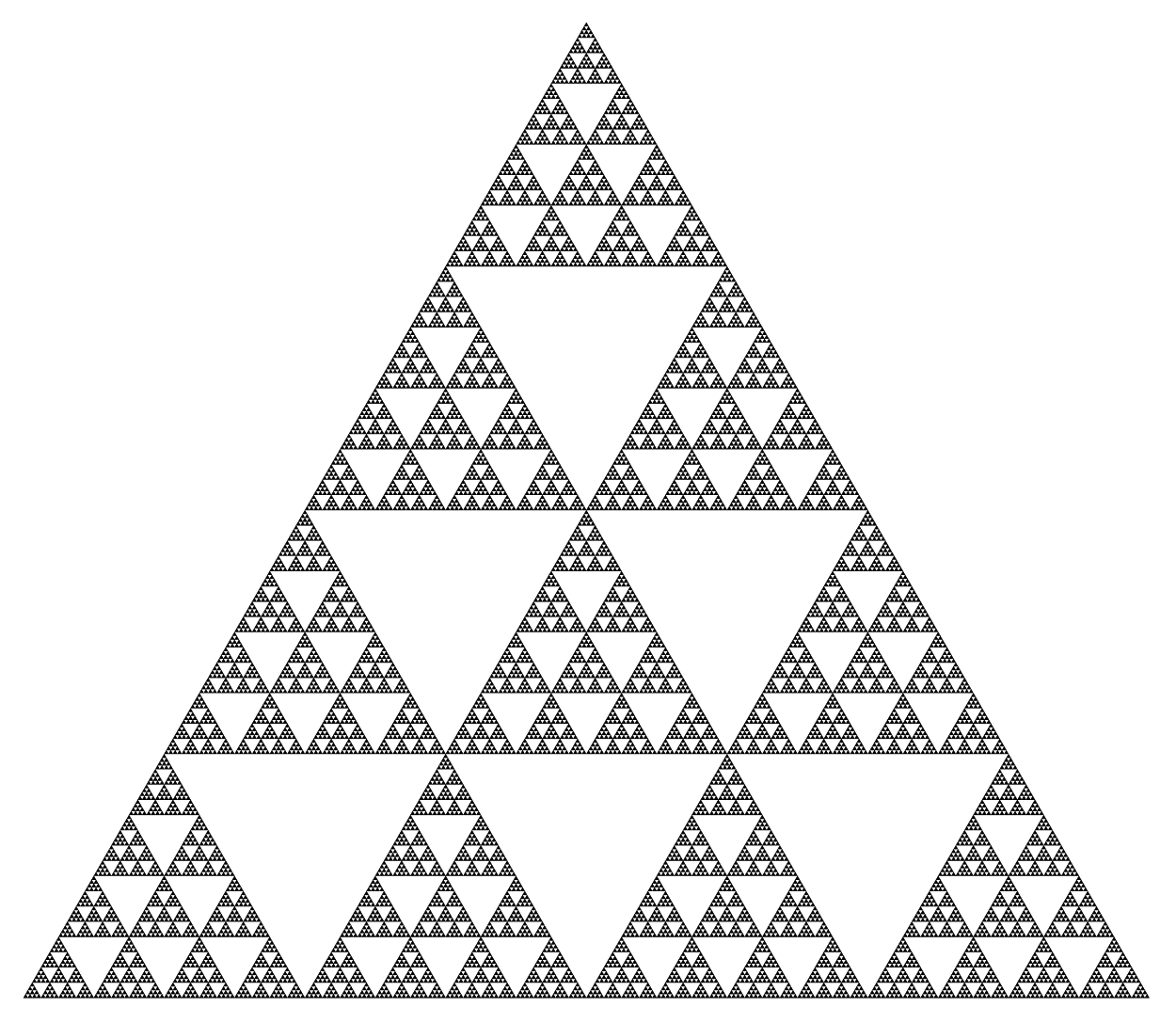}
\setlength{\unitlength}{1cm}
\caption{$\mathcal{SG}_2,\mathcal{SG}_3,\mathcal{SG}_4$.}

\begin{picture}(0,0)
\put(-0.2,4.6){$q_0$}
\put(-2.0,1.1){$q_1$}
\put(1.6,1.1){$q_2$}

\put(-4.5,4.6){$q_0$}
\put(-6.3,1.1){$q_1$}
\put(-2.7,1.1){$q_2$}

\put(4.15,4.6){$q_0$}
\put(2.35,1.1){$q_1$}
\put(5.95,1.1){$q_2$}
\end{picture}
\end{center}
\end{figure}

Recall that $\mathcal{SG}_n$ is the unique nonempty compact subset of $\mathbb{R}^2$ satisfying $\mathcal{SG}_n=\bigcup_{i=0}^{\frac{n^2+n-2}{2}}F_i\mathcal{SG}_n$ with $F_i$'s being contraction mappings defined as $F_i(z)=n^{-1}z+d_{n,i}$ with suitable $d_{n,i}\in \mathbb{R}^2$. See Figure 1.1.  The set $V_0$ consisting of the three vertices $q_0, q_1, q_2$ of the smallest triangle containing $\mathcal{SG}_n$ is called the \textit{boundary}. In this note, we mainly discuss the Laplacians on $V_1=\bigcup_{i=0}^{\frac{n^2+n-2}{2}}V_0$.

Let $D\in \mathcal{LA}(V_0)$, and $\mathcal{E}_0$ be its associated Dirichlet form
$$\mathcal{E}_0(u,v)=\sum_{p,q\in V_0}D_{pq}(u(p)-u(q))(v(p)-v(q)),\quad u,v\in l(V_0).$$
For $\textbf{r}=(r_0,r_1,\cdots,r_{\frac{n^2+n-2}{2}})$ with $r_i>0$, denote by $\mathcal{E}_1$  the induced Dirichlet form on $V_1$ with
\[\mathcal{E}_1(u,v)=\sum_{i=0}^{\frac{n^2+n-2}{2}}r_i^{-1}\mathcal{E}_0(u\circ F_i,v\circ F_i),\quad u,v\in l(V_1).\]
 We call the pair $(D,\textbf{r})$ a \textit{harmonic structure}  if the restriction of $\mathcal{E}_1$ on $V_0$ equals $\mathcal{E}_0$. Once the pair $(D,\textbf{r})$ is a harmonic structure, it could naturally induce a self-similar Dirichlet form on the fractal $\mathcal{SG}_n$, see [K3]. The Laplacian on $V_1$ associated with $\mathcal{E}_1$, denoted by $H_1$, can be expressed as
\[H_1=\sum_{i=0}^{\frac{n^2+n-2}{2}}r_i^{-1} R_i^tDR_i,\]
where $R_i:l(V_1)\rightarrow l(V_0)$ is defined by $R_iv=v\circ F_i$.

By Theorem 3.2.11 in [K3], it is easy to check that $(H_1)_{pq}>0$ if and only if $p,q\in F_iV_0$ for some $0\leq i\leq \frac{n^2+n-2}{2}$ due to the topology of $\mathcal{SG}_n$. We always denote $p\sim q$ if $(H_1)_{pq}>0$. The edge relation '$\sim$' together with the vertex set $V_1$ gives the \textit{level-$1$ graph approximation} $G_1$ of $\mathcal{SG}_n$. See Figure 1.2 for $n=2,3,4$ cases. 

\begin{figure}[h]
\begin{center}
\includegraphics[width=3.7cm]{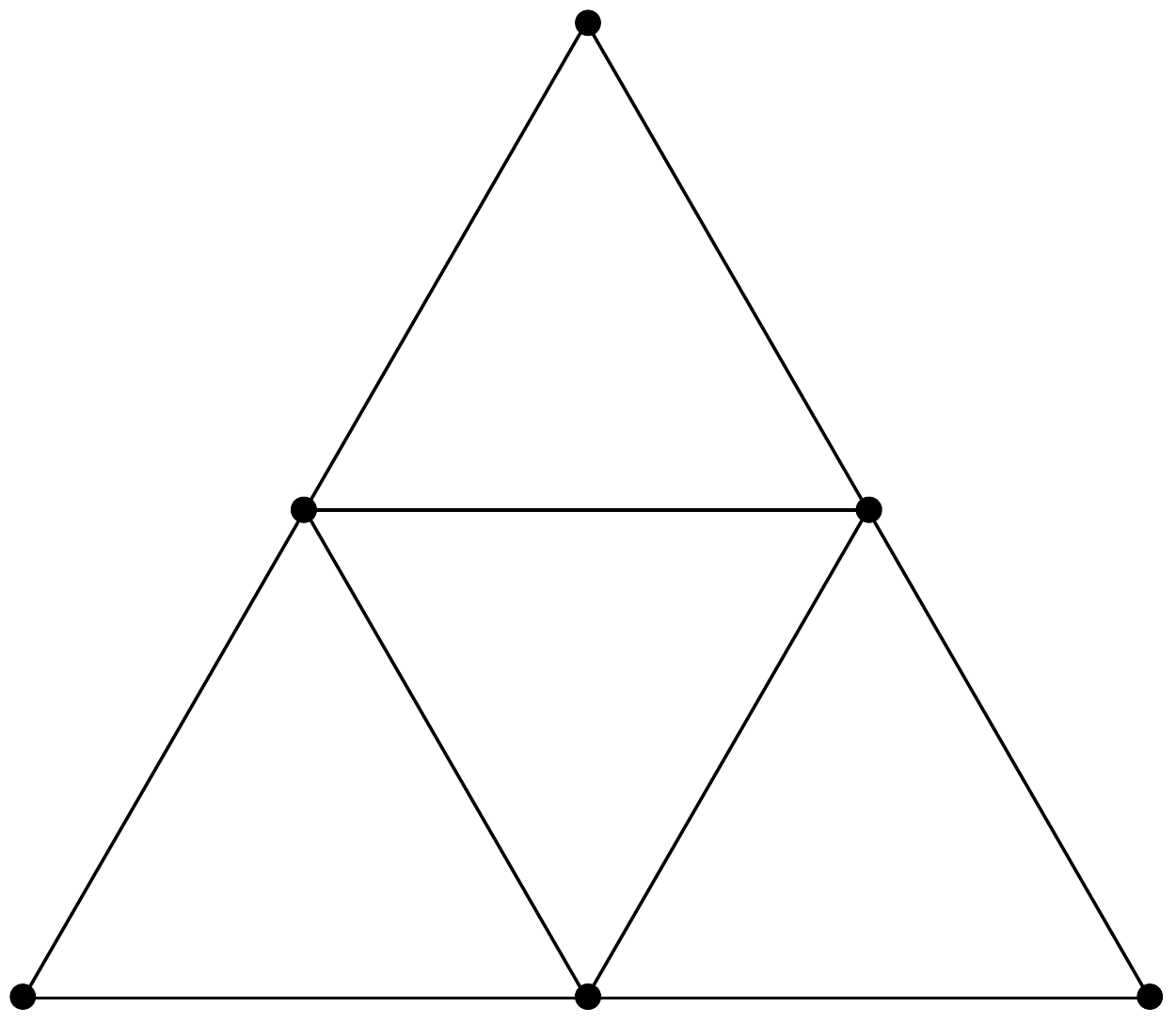}\hspace{0.5cm}
\includegraphics[width=3.7cm]{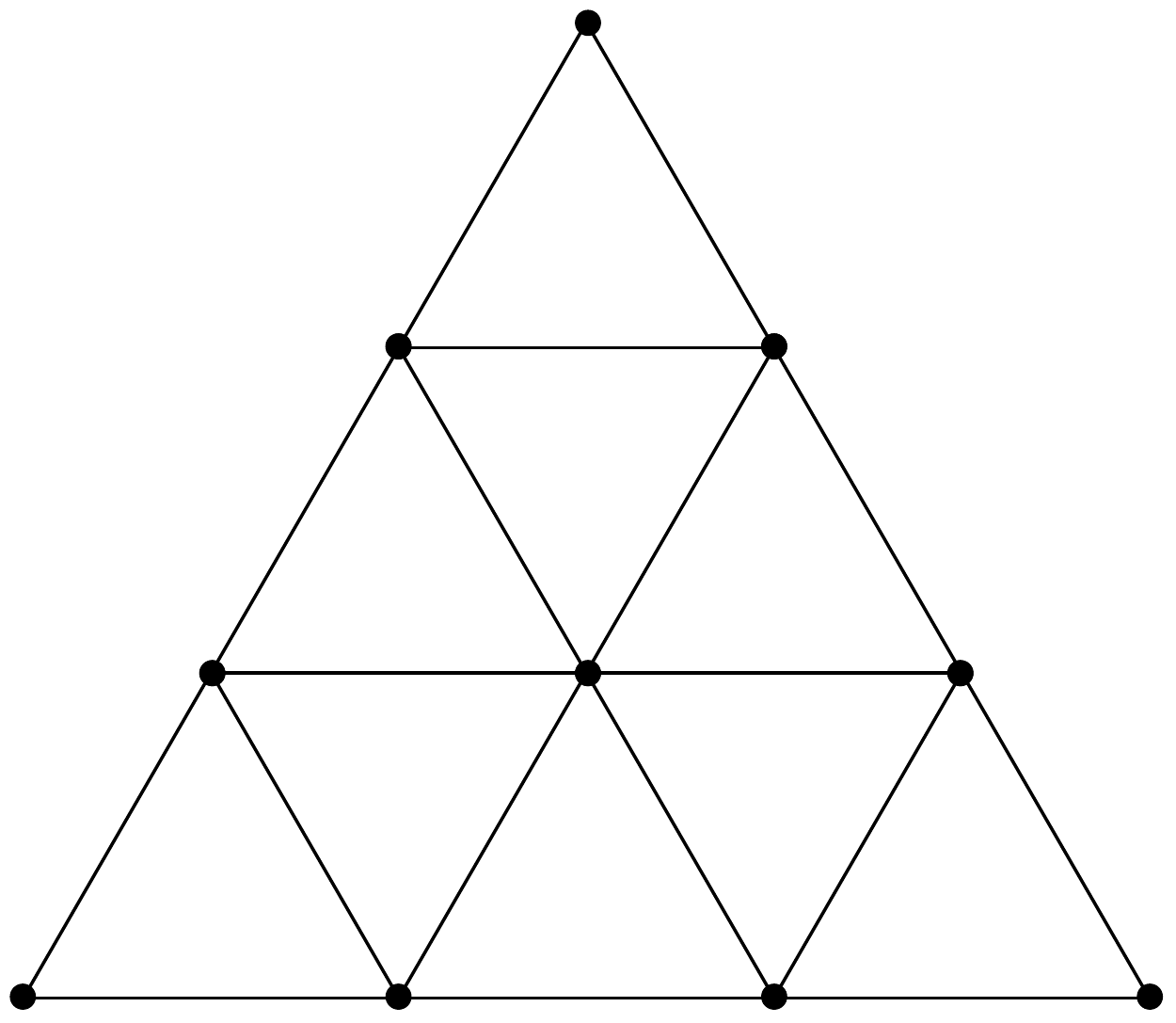}\hspace{0.5cm}
\includegraphics[width=3.7cm]{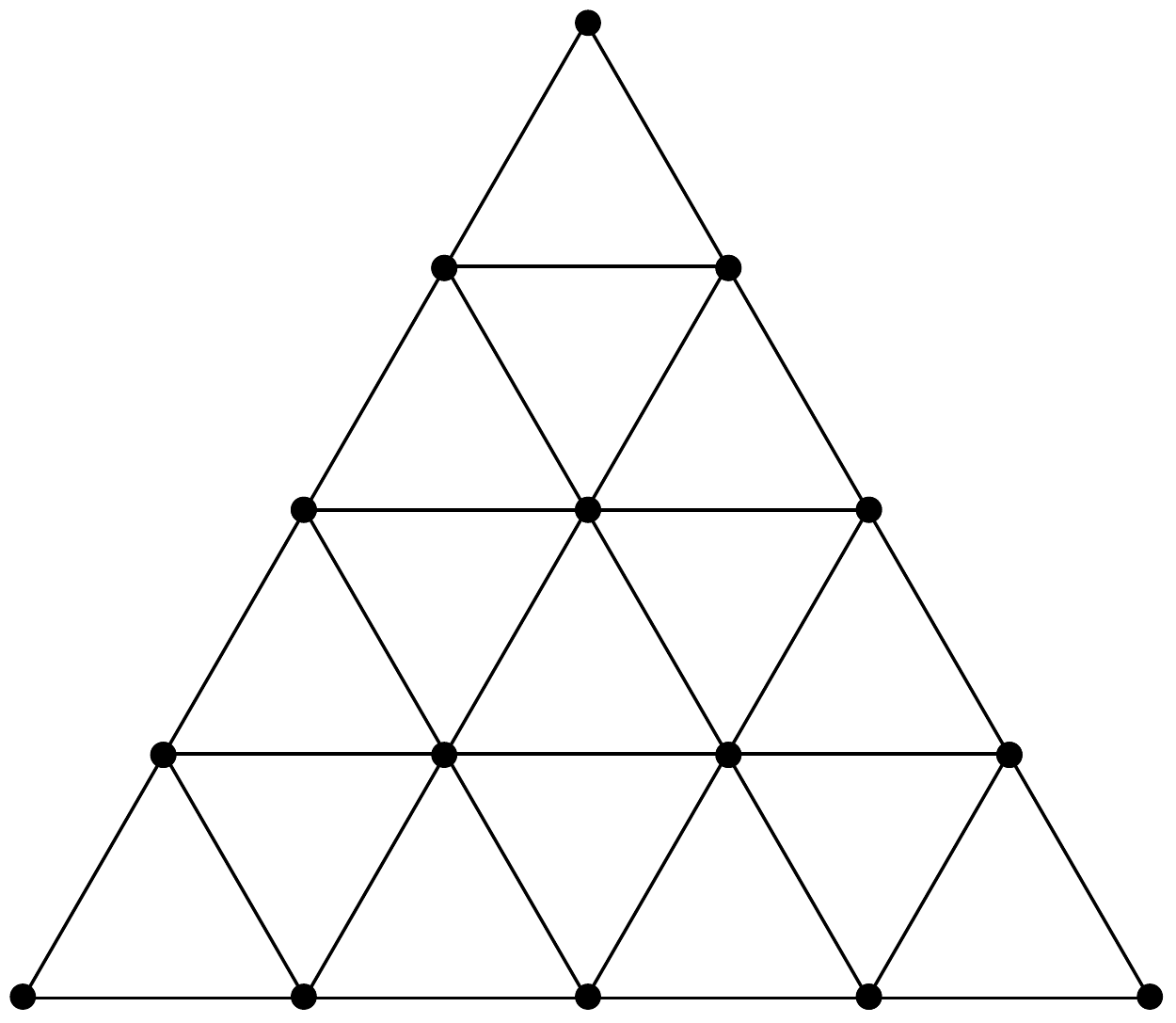}
\setlength{\unitlength}{1cm}
\caption{The level-$1$ graph approximation $G_1$ of $\mathcal{SG}_2,\mathcal{SG}_3,\mathcal{SG}_4$.}

\begin{picture}(0,0)
\put(-0.2,4.6){$q_0$}
\put(-2.0,1.1){$q_1$}
\put(1.6,1.1){$q_2$}

\put(-4.5,4.6){$q_0$}
\put(-6.3,1.1){$q_1$}
\put(-2.7,1.1){$q_2$}

\put(4.15,4.6){$q_0$}
\put(2.35,1.1){$q_1$}
\put(5.95,1.1){$q_2$}
\end{picture}
\end{center}
\end{figure}

Let $h$ be a harmonic function on $\mathcal{SG}_n$ and $v=h|_{V_1}\in l(V_1)$, then 
\begin{equation}
(H_1v)|_{V_1\setminus V_0}=0.
\end{equation}
Conversely, for any $v$ satisfying (1.1), it has a unique extension $h$ so that $h$ is  harmonic on $\mathcal{SG}_n$. Thus we only need to study the solutions of (1.1). The following discrete \textit{maximum principle} is important. See Proposition 2.1.7 in [K3].

\textbf{Proposition 1.3.(Maximum Principle)} \textit{Let $V$ be a finite set and  $H\in \mathcal{LA}(V)$. Let $U$ be a proper subset of $V$. For $p\in V\setminus U$, write
$$\begin{aligned}U_p=\{q\in U:\text{There exist $p_1,p_2,...,p_m\in V\setminus U$ with $p_1=p$}\\ H_{p_ip_{i+1}}>0 \text{ for $i=1,2,...,m-1$ and $H_{p_mq}>0$}\}.\end{aligned}$$
Then if $(Hv)|_{V\setminus U}=0$, then for any $p\in V\setminus U$,
$$\min_{q\in U_p} v(q)\leq v(p)\leq \max_{q\in U_p} v(q).$$
Moreover, $v(p)=\max_{q\in U_p}v(q)$(or $v(p)=\min_{q\in U_p} v(q)$) if and only if $v$ is constant on $U_p$.}

\section{Main Result}
 Call a finite sequence of vertices $\{p_1,p_2,\cdots,p_m\}\subset V_1$ a \textit{chain} if $p_i\sim p_{i+1}$ for $1\leq i<m$ and $p_i\neq p_j$ for $i\neq j$. Before proving Theorem 1.1, we list some lemmas.

\textbf{Lemma 2.1.} \textit{Let $v\in l(V_1)$ satisfying (1.1), and assume $v$ is not a constant. If $p\in V_1\setminus V_0$ is a vertex with a neighbor $q\sim p$ such that $v(q)\neq v(p)$, then there exists a chain $\{p_1,p_2,\cdots,p_m\}$ with $p_1=p$ and $p_m\in V_0$ satisfying
$$v(p_1)<v(p_2)<\cdots<v(p_m).$$
Similarly, there is another chain satisfying $v(p_1)>v(p_2)>\cdots>v(p_m).$
}

\textit{Proof.} Let $p_1=p$. There must exist some $q\sim p_1$ such that $v(p_1)<v(q)$, since by (1.1), $$(H_1v)(p_1)=\sum_{p'\sim p_1} (H_1)_{p_1p'}(v(p')-v(p_1))=0.$$
Set this $q$ as $p_2$. If $p_2\in V_0$, then we have already find the chain. Otherwise, repeat the above procedure until we find a $p_m\in V_0$. The process will stop after finite times of operations, since $V_1$ is a finite set. \hfill$\square$

\textbf{Lemma 2.2.} \textit{Let $A$ be a proper subset of $V_1\setminus V_0$ with at least two vertices. Write $\partial A=\{p\in A|\text{There exists }q\in V_1\setminus (A\bigcup V_0) \text{ such that }(H_1)_{pq}>0\}$. Then $\#\partial A\geq 2$.}

\textit{Proof.} If $\partial A=A$, then we immediately get the lemma. Otherwise, it is easy to see that $A\setminus \partial A$ and $V_1\setminus (A\bigcup V_0)$ are disconnected with each other, which also leads to $\#\partial A\geq 2$ as $(V_1\setminus V_0,\sim)$ is 2-connected. Here we call a graph \textit{$2$-connected} if it could not be disconnected by removing $1$ vertex.   \hfill$\square$

{\textit{Proof of Theorem 1.1.} Let $v$ be the restriction of a nonconstant harmonic function $h$ on $V_1$. For $c\in \mathbb{R}$, define
\[
E(v,c)=\{p\in V_1: v(p)=c\}.\] 
We only need to show that $F_iV_0\nsubseteq E(v,c)$ for any $0\leq i\leq \frac{n^2+n-2}{2}$ and any $c\in \mathbb{R}$.

Assume $F_{i_1}V_0\subset E(v,c)$ for some $i_1$ and $c$. We may find another $i_2$ such that $F_{i_2}V_0\subset E(v,c)$ and $F_{i_2}V_0\bigcap F_{i_1}V_0\neq \emptyset$. Also, another $i_k$ such that $F_{i_k}V_0\subset E(v,c)$ and $F_{i_k}V_0\bigcap (\bigcup_{l=1}^{k-1} F_{i_l}V_0)\neq \emptyset$. Repeat the procedure until we can not find another such $i$, and set $A$ to be the union of these $F_{i_k}V_0$'s. As a consequence of Lemma 2.2, there should be at least two vertices, denoted by $a_1,a_2\in A\setminus V_0$, such that for each $a_i$, there exists a neighbor $p\sim a_i$ with $v(p)\neq c$.

Let $CH=\{p_1,p_2,\cdots,p_m\}\subset A$ be a \textit{geodesic chain} (with minimal number of edges, may not be unique) in $A$ connecting $p_1=a_1$ and $p_m=a_2$. See Figure 2.1 for examples of such chains. The value of $v$ is always $c$ along this chain. On the other hand, for each $i=1,2$, by using Lemma 2.1, we could connect $a_i$ to the boundary $V_0$ by two different chains $CH^{S,i}$ and $CH^{L,i}$, denoted as
$$\begin{aligned}
CH^{S,i}=\{p^{S,i}_k\}_{k=1}^{m_{S,i}},\text{ with }p_1^{S,i}=a_i\text{ and }p^{S,i}_{m_{S,i}}\in V_0,\\
CH^{L,i}=\{p_k^{L,i}\}_{k=1}^{m_{L,i}},\text{ with }p_1^{L,i}=a_i\text{ and }p^{L,i}_{m_{L,i}}\in V_0,\end{aligned}
$$
satisfying  that $p_k^{S,i}<c,\forall 2\leq k\leq m_{S,i}$ and $p_k^{L,i}>c,\forall 2\leq k\leq m_{L,i}.$

\begin{figure}[h]
\begin{center}
\includegraphics[width=4cm]{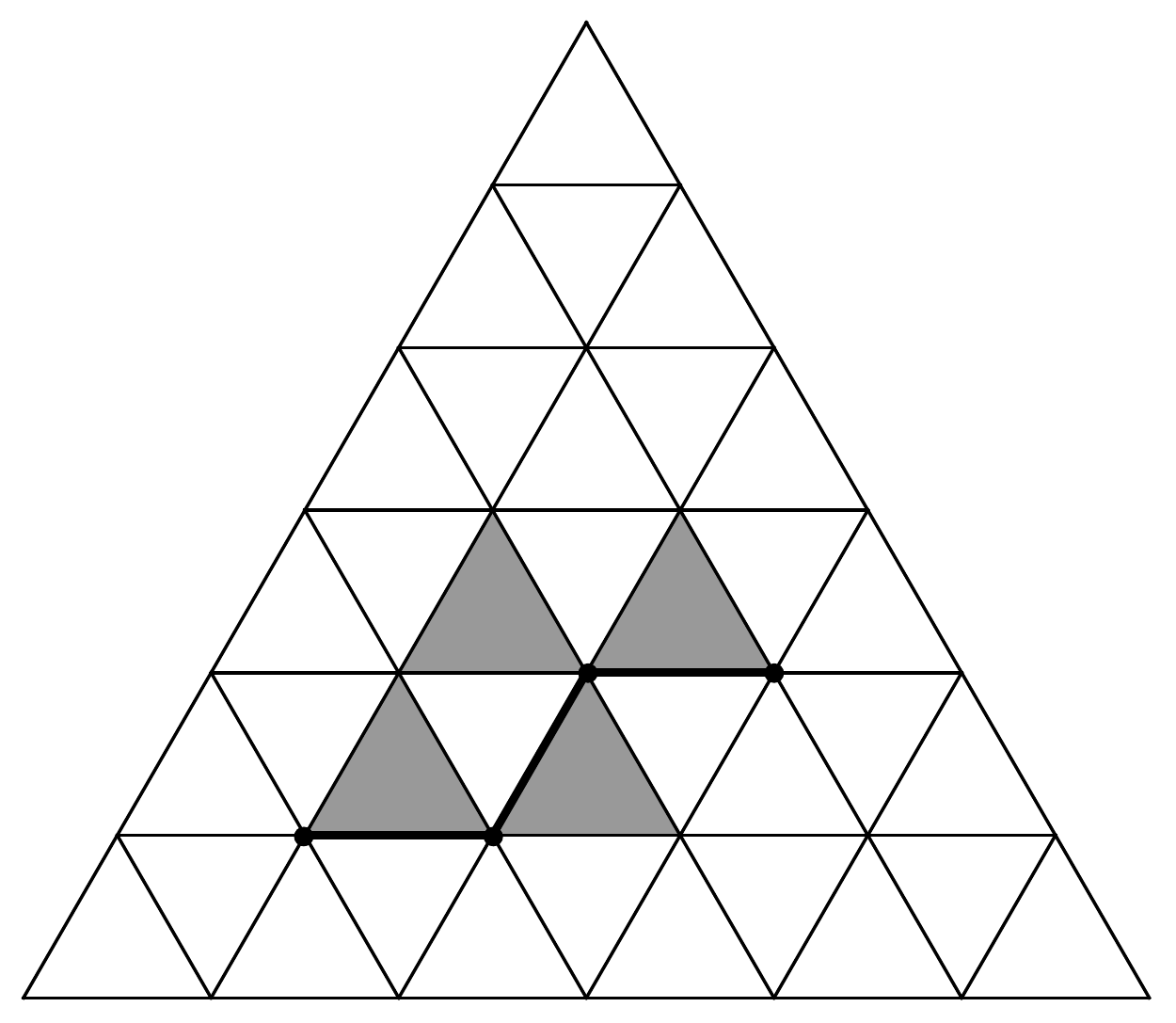}
\includegraphics[width=4cm]{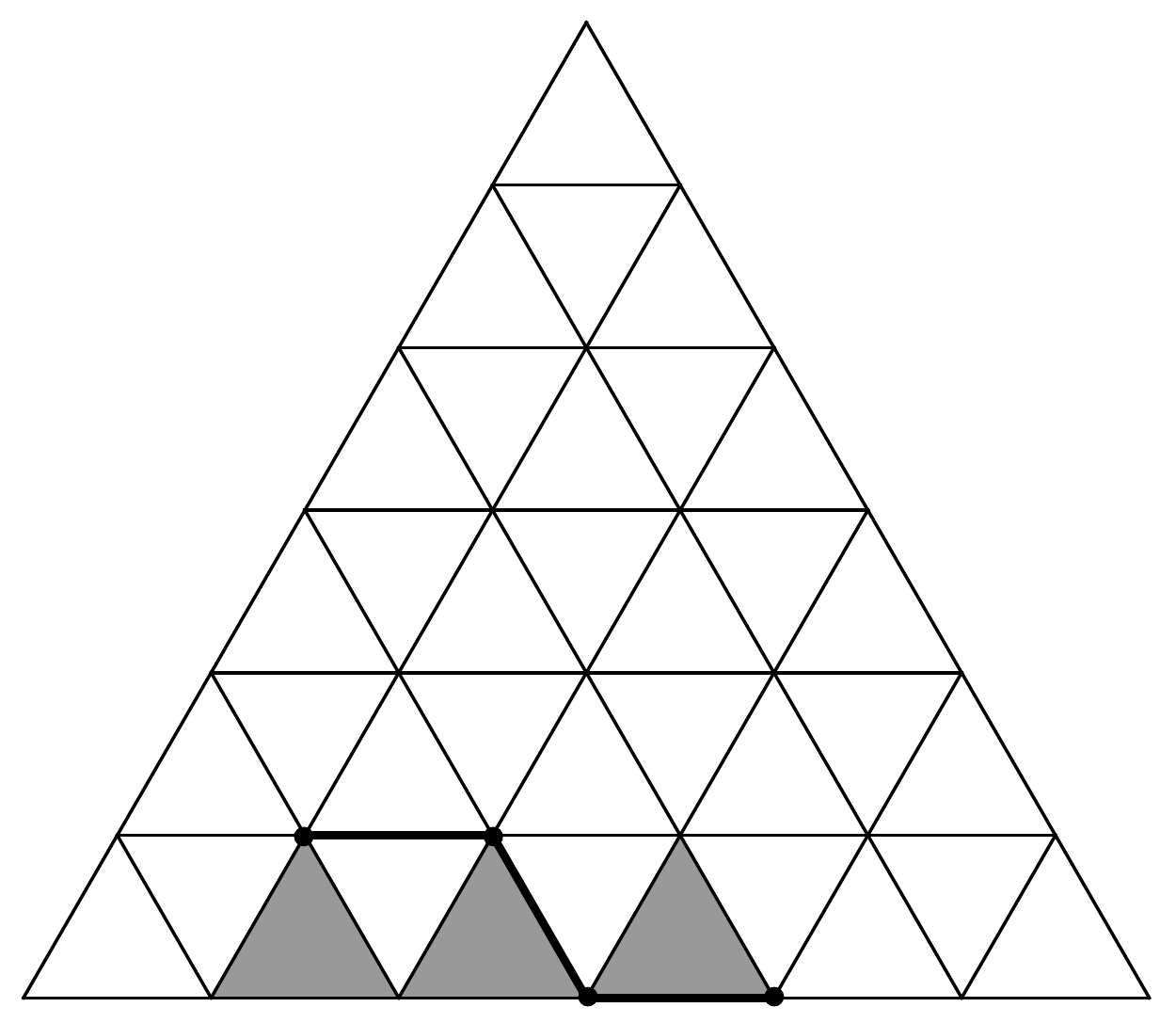}
\includegraphics[width=4cm]{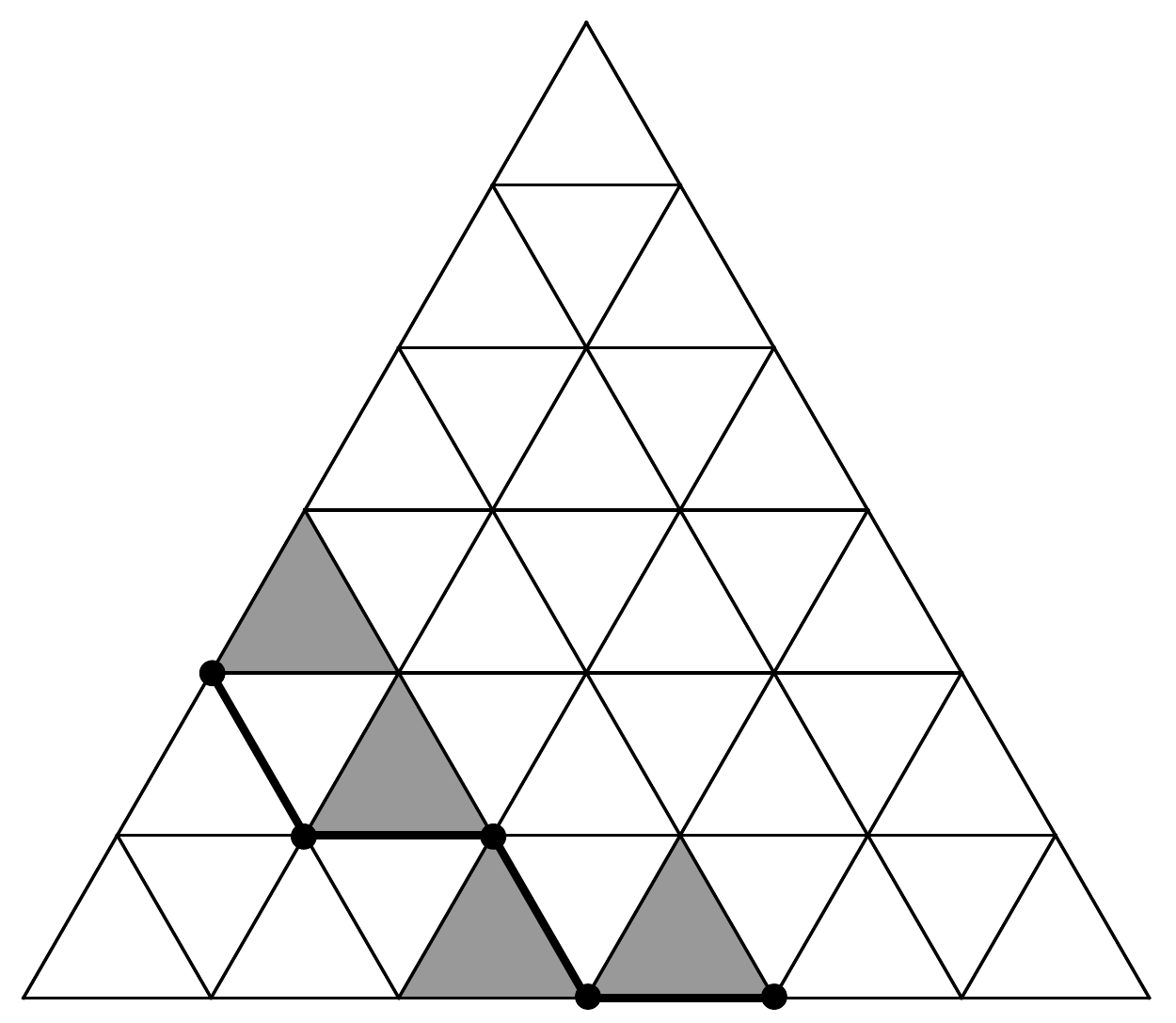}
\setlength{\unitlength}{1cm}
\caption{Examples of possible cases for $CH$.}
\end{center}
\end{figure}

Take a vertex $p\in A\setminus CH$ such that $p\sim p_{i_1},p\sim p_{i_2}$ for two different $p_{i_1},p_{i_2}\in CH$. Note that $v(p)=c$. It is not possible for $p\in V_0$, otherwise it will  contradict to the fact that there exist five chains as shown above, since $CH^{S,i}$ and $CH^{L,j}$ never cross when $i\neq j$. Now as in Proposition 1.3, we take $U$ to be the union of $V_0$ and above five chains, and $V$ to be $V_1$. Noticing that $p\in V\setminus U$, we have
\begin{equation}
\text{ either }U_p\subset CH\cup CH^{S,1}\cup CH^{S,2}\text { or } U_p\subset CH\cup CH^{L,1}\cup CH^{L,2}.
\end{equation}
Moreover, since $CH$ is a geodesic chain, it could not encircle a nonempty subset of $V_1\setminus V_0$,  and thus there must exist a chain in $V_1\setminus CH$ which connects $p$ to some $q_i\in V_0$.  So we have
\begin{equation}
U_p\setminus CH\neq\emptyset.
\end{equation}

By (2.1), (2.2) and the definition of $CH^{S,i}$ and $CH^{L,i}$, we have either $\min_{q\in U_p} v(q)<c=\max_{q\in U_p}v(q)$ or $\min_{q\in U_p}v(q)=c<\max_{q\in U_p}v(q)$. Thus we get $v(p)\neq c$ by using the maximum principle as stated in Proposition 1.3. However, this contradicts to the fact that $p\in A$.  $\hfill\square$

 As we know, harmonic functions on $\mathcal{SG}_n$ are of $3$-dimensional, and could be uniquely determined by their boundary values on ${V_0}$. More precisely, there exist \textit{harmonic extension matrices} $A_i$, $0\leq i\leq\frac{n^2+n-2}{2}$, so that \[h|_{F_iV_0}=A_ih|_{V_0}\]
 holds for any harmonic function $h$ on $\mathcal{SG}_n$.
We could restate Theorem 1.1 as following.  

\textbf{Corollary 2.3.}  \textit{All the harmonic extension matrices $A_i$ of $\mathcal{SG}_n$ are nondegenerate.}

\textit{Proof.} Let $h$ be a harmonic function and let $v=h|_{V_1}$. If $v|_{V_0}=c\neq 0$, then $v|_{F_iV_0}=c$. If $v|_{V_0}$ is not a constant, then by Theorem 1.1,  $v|_{F_iV_0}$ is also not a constant. These show that $A_iv|_{V_0}\neq 0$ if $v|_{V_0}\neq 0$. $\hfill\square$

\end{document}